\documentclass[a4paper, 12pt]{article}

\usepackage{amsfonts, amsmath, amsthm}
\usepackage{latexsym}

\theoremstyle{plain}
\newtheorem{thm}{Theorem}[section]
\newtheorem{prp}{Proposition}[section]
\newtheorem{lmm}[prp]{Lemma} 
\newtheorem{crl}[prp]{Corollary} 

\theoremstyle{definition}

\theoremstyle{remark}
\newtheorem{rem}{Remark}[section]

\numberwithin{equation}{section}

\usepackage{amssymb}
\usepackage{latexsym}

\newcommand{\R}{\mathbb{R}}
\newcommand{\C}{\mathbb{C}}

\newcommand{\pa}{\partial}
\newcommand{\eps}{\varepsilon}
\newcommand{\diag}{\rm diag}

\newcommand{\dirac}{\mathcal{D}}
\newcommand{\conjdirac}{\widetilde{\mathcal{D}}}
\newcommand{\OOmega}{\mbox{\boldmath$\Omega$}}

\begin{document}
\title{
A remark on the algebraic normal form method \\
applied to the Dirac-Klein-Gordon system \\
in two space dimensions
}  

\author{
         Masahiro Ikeda\footnote{
          Department of Mathematics, Osaka University, 
          Toyonaka, Osaka 560-0043, Japan.
       \newline
         E-mail: \texttt{m-ikeda@cr.math.sci.osaka-u.ac.jp}}
   \and  
         Akihiro Shimomura\footnote{  
        Graduate School of Mathematical Sciences, The University of Tokyo, 
        3-8-1 Komaba, Meguro-ku, Tokyo, 153-8914 Japan. 
        E-mail:\ {\tt simomura@ms.u-tokyo.ac.jp}} 
   \and  
         Hideaki Sunagawa\thanks{
             Department of mathematics, Osaka University, 
             Toyonaka, Osaka 560-0043, Japan. 
        \newline
            E-mail: \ {\tt sunagawa@math.sci.osaka-u.ac.jp}}
 } 
 
\date{\today}   
\maketitle

\noindent{\bf Abstract:}\ 
We consider the massive Dirac-Klein-Gordon system in two space dimensions. 
 Under the non-resonace mass condition, we show that the solution is 
 asymptotically free if the initial data are sufficiently small in a suitable 
 weighted Sobolev space. In particular, it turns out that the Dirac component 
 of the DKG system tends to a solution of the free Dirac 
 equation. Our proof is based on the algebraic normal form method. 
\\



\section{Introduction} 

This paper is intended to give a remark on applications of the algebraic 
normal form method developed by \cite{shatah}, \cite{kosecki}, \cite{o_t_t}, 
\cite{katayama}, \cite{tsutsumi}, \cite{sunagawa}, \cite{k_o_s}, etc. 
The model equation which we focus on is the two-dimensional massive 
Dirac-Klein-Gordon system 
\begin{align}
\left\{
\begin{array}{l}
\dirac_M \psi=ig \phi \beta \psi,\\
(\Box+m^2)\phi=g \langle \psi, \beta \psi\rangle_{\C^2},
\end{array}\right.
\quad (t,x)\in \R\times \R^2
 \label{DKG}
\end{align}
with the initial condition 
\begin{align}
 (\psi, \phi, \pa_t\phi)|_{t=0} = (\psi_0, \phi_0,\phi_1), \qquad x \in \R^2.
 \label{data_DKG}
\end{align}
Here $(\psi,\phi)$ is a $\C^2\times \R$-valued unknown function of 
$(t,x) \in \R\times \R^2$. $M$, $m$ are positive constants, $g$ is a real 
constant, $\Box=\pa_t^2-\Delta$, $\Delta=\pa_1^2+\pa_2^2$, $\pa_j=\pa/\pa x_j$ 
($j=1,2$) and $\langle \,\cdot ,\cdot\, \rangle_{\C^2}$ denotes the 
standard scalar product in $\C^2$, i.e., 
$\langle u,v\rangle_{\C^2}=u^{\dagger} v$ for $u, v \in \C^2$ 
(regarded as column vectors), 
where $u^{\dagger}$ is  the complex conjugate transpose of $u$. 
The Dirac operator $\dirac_M$ is defined by 
$$
 \dirac_{M}
 =\pa_t+\alpha_1\pa_1+\alpha_2 \pa_2 +iM\beta
 =\pa_t+\alpha \cdot \nabla_x +iM\beta
$$ 
with $2\times 2$ hermitian matrices $\alpha_1$, $\alpha_2$, $\beta$ satisfying 
\begin{align*}
 &\alpha_1^2=\alpha_2^2=\beta^2=I
 ,\\
 &\alpha_1\alpha_2+\alpha_2\alpha_1=
 \alpha_1\beta+\beta \alpha_1=\alpha_2\beta+\beta \alpha_2=O
 .
\end{align*}
We also set 
$$
 \conjdirac_{M}
 =\pa_t-(\alpha \cdot \nabla_x +iM\beta),
$$ 
then we can easily check that the following relations hold: 
$$
 \dirac_M\conjdirac_M=\conjdirac_M\dirac_M=(\Box+M^2)I.
$$
This implies that the solution $(\psi,\phi)$ of (\ref{DKG})--(\ref{data_DKG}) 
also solves 
\begin{align}
\left\{
\begin{array}{l}
(\Box+M^2) \psi=ig \conjdirac_M (\phi \beta \psi),\\
(\Box+m^2)\phi=g \langle \psi, \beta \psi \rangle_{\C^2},
\end{array}\right.
\quad (t,x)\in \R\times \R^2
 \label{NLKG}
\end{align}
with the initial condition 
\begin{align}
 (\psi, \pa_t \psi, \phi, \pa_t\phi)|_{t=0} = 
 (\psi_0, \psi_1,\phi_0,\phi_1), \qquad x \in \R^2,
 \label{data_NLKG}
\end{align}
where 
$\psi_1
=-(\alpha \cdot \nabla_x + iM\beta) \psi_0 +ig \phi_0\beta \psi_0$. 
According to Theorem~6.1 of \cite{sunagawa} (see also \cite{tsutsumi}, 
\cite{k_o_s}), the solution of (\ref{NLKG})-(\ref{data_NLKG}) exists globally 
in time if $m\ne 2M$ and the data are sufficiently small, smooth and 
decay fast as $|x|\to \infty$. 
Moreover there exists a solution ($\psi^{\pm}, \phi^{\pm}$) of the 
free Klein-Gorodon equation 
$$
\left\{
\begin{array}{l}
(\Box+M^2) \psi^{\pm}=0,\\
(\Box+m^2)\phi^{\pm}=0,
\end{array}\right.
$$
such that 
$$
 \lim_{t\to \pm \infty} 
 \sum_{j=0}^{1}
 \Bigl(\|\pa_t^j(\psi(t,\cdot)-\psi^{\pm}(t,\cdot))\|_{H^{1-j}}
 +
 \|\pa_t^j(\phi(t,\cdot)-\phi^{\pm}(t,\cdot))\|_{H^{1-j}}
 \Bigr)
 =0.
$$
In this sense, the solution of (\ref{DKG})--(\ref{data_DKG}) behaves like 
a solution of the free Klein-Gordon equations in the large time if $m\ne 2M$. 
However, this does {\em not} imply that the solution is asymptotically free.  
What we emphasize here is that a solution $u$ of the free 
Klein-Gordon equation $(\Box+M^2) u=0$ is not necessarily a solution 
of the free Dirac equation $\dirac_M u=0$ in general. 
So the following question arises: 
{\em Does the Dirac componet $\psi(t)$ of (\ref{DKG}) tend to a solution 
of the free Dirac quation as $t \to \pm \infty$?} 
As far as the authors know, there are no previous papers which address 
this question in the case of two space dimensions. 
There are several results in 3D case (see e.g., \cite{h_i_n} and 
the referances therein), however, those methods do not work well in 2D case 
because of the insufficiency of expected decay rate with respect to $t$ of 
the nonlinear terms. 
We will give an affirmative answer to this question by using the algebraic 
normal form method.

To state the main result, let us introduce the weighted Sobolev space 
$$
 H^{s,k}(\R^2)
 =
 \{ 
  u \in L^2(\R^2)\, :\, 
  (1+|\cdot|^2)^{k/2} (1-\Delta)^{s/2} u \in L^2(\R^2) 
 \}
$$
equipped with the norm 
$$
 \|u\|_{H^{s,k}(\R^2)}
 = \|(1+|\cdot|^2)^{k/2} (1-\Delta)^{s/2} u\|_{L^2(\R^2)}.
$$
As usual, we write $H^s=H^{s,0}$ and $\|u\|_{H^s}=\|u\|_{H^{s,0}}$. 
Our main result is as follows.
\begin{thm}\label{thm_main}
Let $m\ne 2M$. Assume that  
$(\psi_0, \phi_0,\phi_1) 
\in H^{s+1,s}\times H^{s+1,s}\times H^{s,s}(\R^2)$ 
with $s \ge 18$. 
There exists a positive constant $\eps$ such that if 
\begin{equation}
\label{DataSize}
 \|\psi_0\|_{H^{s+1,s}(\R^2)}+ \|\phi_0\|_{H^{s+1,s}(\R^2)}+
 \|\phi_1\|_{H^{s,s}(\R^2)} \le \eps,
\end{equation}
the Cauchy problem \eqref{DKG}--\eqref{data_DKG} 
admits a unique global solution 
$(\psi,\phi)$ satisfying 
$$
 \psi\in C(\R;H^{s+1}(\R^2)), \quad 
 \phi\in \bigcap_{k=0}^{1}C^{k}([0,\infty);H^{s+1-k}(\R^2)).
$$
Furthermore, there exist 
$\psi_0^{\pm} \in H^{s-1}(\R^2)$ and 
$(\phi_0^{\pm},\phi_1^{\pm}) \in H^{s-1}\times H^{s-2}(\R^2)$
such that 
$$
 \lim_{t\to \pm \infty} 
 \|\psi(t,\cdot)-\psi^{\pm}(t,\cdot)\|_{H^{s-1}}
 =0
\quad \mbox{and}\quad   
 \lim_{t\to \pm \infty} 
 \sum_{j=0}^{1}\|\pa_t^j(\phi(t,\cdot)-\phi^{\pm}(t,\cdot))\|_{H^{s-1-j}}
 =0,
$$
where $\psi^{\pm}$ and $\phi^{\pm}$ are the solutions to  
$$ 
 \left\{\begin{array}{l}
 \dirac_M \psi^{\pm}=0\\
 \psi^{\pm}|_{t=0}=\psi_0^{\pm}
 \end{array}\right.
\quad \mbox{and}\quad   
 \left\{\begin{array}{l}
 (\Box+m^2)\phi^{\pm}=0\\
 (\phi^{\pm}, \pa_t \phi^{\pm})|_{t=0}=(\phi_0^{\pm}, \phi_1^{\pm}),
 \end{array}\right.
$$ 
respectively.
\end{thm}

\begin{rem}
The condition $m\ne 2M$ is often called the non-resonance mass condition. 
Difficulties appearing in the resonant case ($m=2M$) are explained in 
\cite{kawasuna}.
\end{rem}

\begin{rem}
Recently, the first author considered the final state problem for \eqref{DKG} 
in two space dimensions and succeeded in showing the existence of wave 
operators for \eqref{DKG} under the non-resonance mass condition. 
See \cite{ikeda} for the detail.
\end{rem}

The rest of this paper is organized as follows: 
In the next section, we give some preliminaries mainly on the commuting vector 
fields and the null forms. In Section~3, we recall and develop an algebraic 
normal form transformation. We will get an a priori esimate of the solution in 
Section~4. After that, Theorem \ref{thm_main} will be proved in Section~5. 
Throughout this paper, we will frequently use the following 
conventions on implicit constants: 
\begin{itemize}
\item
$A\lesssim B$ (resp. $A \gtrsim B$) stands for  $A \le CB$ (resp. $A \ge CB$) 
with a positive constant $C$. 
\item 
The expression $f=\sum_{\kappa \in K}' g_\kappa$ means that there exists a 
family $\{C_\kappa\}_{\kappa \in K}$ of constants such that 
$f=\sum_{\kappa \in K} C_\kappa g_\kappa$. 
\end{itemize}
Also, the notation $\langle y \rangle=(1+|y|^2)^{1/2}$ will be used for 
$y \in \R^N$ with a positive integer $N$.

\section{Commuting vector fields and the null forms} 
In this section, we summarize basic properties of some vector fields 
associated with the Klein-Gordon operators.  
We put
$x_0=-t$, $x=(x_1,x_2)$,  
$\Omega_{ab} = x_{a}\pa_{b}-x_{b}\pa_{a}$, $0\le a, b \le 2$, 
$\pa=(\pa_0,\pa_1,\pa_2)=(\pa_t,\pa_{x_1},\pa_{x_2})$ and 
$$
 Z=(Z_1,\dots, Z_{6})
 =\bigl( \pa_0,\pa_1,\pa_2,\Omega_{01}, \Omega_{02}, \Omega_{12} \bigr).
$$
Note that the following commutation relations hold: 
\begin{align}
\label{commutation}
&[\Box +m^2, Z_j] = 0, 
\\
&[\Omega_{a b}, \pa_{c}]  = \eta_{b c} \pa_{a} - \eta_{c a} \pa_{b}, 
\nonumber\\
&[\Omega_{a b}, \Omega_{c d}] 
  =  \eta_{a d} \Omega_{b c}+ \eta_{b c} \Omega_{a d}
   - \eta_{a c} \Omega_{b d} - \eta_{b d} \Omega_{a c}
\nonumber
\end{align}
for $m \in \R$, $1\leq j\leq 6$, $0\leq a, b \leq 2$. 
Here 
$[\cdot ,\cdot ]$ denotes the commutator of linear operators,
and 
$(\eta_{a b})_{0\leq a,b \leq 2}= \diag (-1,1,1)$. 
Note that $\Box=-\sum_{a,b=0}^{2}\eta_{ab}\pa_a\pa_b$. 
For a smooth function $u$ of $(t,x) \in \R^{1+2}$ and for a non-negative 
integer $s$, we define 
$$
|u(t,x)|_{s} := 
   \sum_{|\nu|\leq s}|Z^{\nu}u(t,x)|
$$
and 
$$
\|u(t)\|_{s} := 
   \sum_{|\nu|\leq s}\|Z^{\nu} u(t,\cdot)\|_{L^2(\R^2)}, 
$$
where $\nu=(\nu_1,\dots, \nu_{6})$ is a multi-index, 
$Z^{\nu} =Z_1^{\nu_1}\cdots Z_{6}^{\nu_{6}}$ and 
$|\nu|= \nu_{1}+\cdots +\nu_{6}$. 
Next we introduce the null form $Q_0$ and the strong null forms $Q_{ab}$ 
as follows:
\begin{align}
 &Q_0(u,v)
 =-\sum_{a,b=0}^{2} \eta_{ab}(\pa_a u)(\pa_b v),
 \label{null}\\
  &Q_{a b}(u,v)
 =(\pa_{a} u) (\pa_{b} v)-(\pa_{b} u) (\pa_{a} v),
 \quad  0\le a,b\le 2.
  \label{strongnull} 
\end{align}
We summarize well known properties on the strong null forms. 

\begin{lmm} \label{lem_null}
Let $u$, $v$ be smooth functions of $(t,x) \in \R^{1+2}$. We have 
\begin{align*}
 |Q_{a b}(u,v)|  \lesssim \frac{1}{\langle |t|+|x| \rangle} 
   \bigl( |u|_1 |\pa v| + |\pa u| |v|_1 \bigr) 
\end{align*}
for $0\le a,b\le 2$, and
\begin{align*}
 Z^{\nu} Q_{a b} (u,v) = 
  \sum_{c,d =0}^{2} \mathop{{\;\,\sum}'}_{|\lambda|+|\mu|\leq |\nu|}
  Q_{c d} (Z^{\lambda} u, Z^{\mu} v) 
\end{align*}
for any multi-index $\nu$.
\end{lmm}

\section{Algebraic normal form transformation} 
This section is devoted to some decomposition of the nonlinear terms in 
(\ref{DKG}). 

Let $v_j$ and $\widetilde{v}_j$ be smooth functions of 
$(t,x) \in \R^{1+2}$ (not necessarily scalar-valued), and let $m_1$, $m_2$ be real constants. We set 
$h_j=(\Box+m_j^2)v_j$ and $\widetilde{h}_j=(\Box+m_j^2)\widetilde{v}_j$ 
for $j=1,2$. 
We write 
$$F\sim G$$ 
if $F-G$ can be written as a linear combination of 
$Q_{ab}(\pa^\mu v_k, \pa^\nu \widetilde{v}_l)$, 
$(\pa^\mu v_k) (\pa^{\nu} \widetilde{h}_l)$, 
$(\pa^\mu h_k)(\pa^\nu \widetilde{v}_l)$ or $h_k\widetilde{h}_l$ 
with $|\mu|$, $|\nu| \le 1$, $0\le a,b \le 2$ and $1\le k, l \le 2$. 
The following lemma is important for our main purpose.

\begin{lmm}[\cite{sunagawa}, \cite{k_o_s}] \label{lem_algebra}
Put $\mathbf{e}_{kl}=v_k \widetilde{v}_l$, 
$\mathbf{\tilde{e}}_{kl}=Q_0(v_k, \widetilde{v}_l)$ 
and $\mathcal{L}_j=\Box+m_j^2$, where $Q_0$ is given by (\ref{null}). 
We have 
$$
 \bigl( 
 \mathcal{L}_j(\mathbf{e}_{kl})\ \ \mathcal{L}_j(\mathbf{\tilde{e}}_{kl})
 \bigr) 
 \sim 
 (\mathbf{e}_{kl}\ \ \mathbf{\tilde{e}}_{kl}) A_{jkl},
$$
where 
\begin{align*} 
A_{jkl}=\begin{pmatrix}
 m_j^2-m_k^2-m_l^2 & 2m_k^2m_l^2\\
 2 & m_j^2-m_k^2-m_l^2 
 \end{pmatrix}.
\end{align*}
\end{lmm}

See Proposition 4.1 of \cite{k_o_s} or Lemma 6.1 of \cite{sunagawa} 
for the proof of this lemma. 
Remark that the proof remains valid in the vector-valued case.

Now we focus our attention to the structure of the matrix $A_{jkl}$. 
Since 
$$
 \det A_{jkl}
 =
 \prod_{\sigma_1,\sigma_2 \in \{\pm 1\}} (m_j+\sigma_1 m_k+\sigma_2 m_l),
$$
we see that $A_{121}$ and $A_{211}$ are invertible if $m_2\ne 2m_1$. 
Moreover we have 
\begin{align*}
 v_k \widetilde{v}_l
 &=
 (\mathbf{e}_{kl}\ \ \mathbf{\tilde{e}}_{kl})
 \begin{pmatrix}1 \\ 0 \end{pmatrix} 
 \nonumber\\
 &=
 (\mathbf{e}_{kl}\ \ \mathbf{\tilde{e}}_{kl})
 A_{jkl}\begin{pmatrix}p_{jkl} \\ \tilde{p}_{jkl} \end{pmatrix}
 \nonumber\\
 &\sim
  (\mathcal{L}_j(\mathbf{e}_{kl})\ \ \mathcal{L}_j(\mathbf{\tilde{e}}_{kl}))
 \begin{pmatrix}p_{jkl} \\ \tilde{p}_{jkl} \end{pmatrix}
 \nonumber\\
  &=
  (\Box+m_j^2)\Bigl(p_{jkl} v_k\widetilde{v}_l 
  + \tilde{p}_{jkl} Q_0(v_k, \widetilde{v}_l) \Bigr)
\end{align*}
with 
$$
 \begin{pmatrix} p_{jkl}\\ \tilde{p}_{jkl} \end{pmatrix}
 =A_{jkl}^{-1} \begin{pmatrix} 1\\ 0 \end{pmatrix}
$$
for $(j,k,l)=(1,2,1)$ or $(2,1,1)$. By using the above formula with 
$(m_1,m_2,v_{2},\widetilde{v}_{1})=(M,m,\phi,\beta \psi)$ 
or $(m_1,m_2,v_{1},\widetilde{v}_{1})=(M,m,\psi^{\dagger},\beta \psi)$, 
we arrive at the following decompositions for the nonlinear terms 
in (\ref{DKG}) : 

\begin{crl} \label{cor_normalform}
Let $(\psi,\phi)$ be a solution for (\ref{DKG}) with $m\ne 2M$. We have 
$$
\left\{\begin{array}{l} 
ig \phi \beta \psi
=\dirac_M(\conjdirac_M\Lambda_D) +N_D +R_D,\\
g \langle \psi, \beta \psi \rangle_{\C^2}
 =(\Box+m^2)\Lambda_{KG} +N_{KG} +R_{KG},
\end{array}\right.
$$
where
\begin{align*} 
 &\Lambda_D=
  \mathop{{\;\,\sum}'}_{|\mu|,|\nu|\le 1}
  (\pa^{\mu}\phi) \beta \pa^{\nu} \psi,\\
 &\Lambda_{KG}=
  \mathop{{\;\,\sum}'}_{|\mu|,|\nu|\le 1}
  \langle \pa^{\mu}\psi, \beta \pa^{\nu} \psi\rangle_{\C^2},\\
 &N_D=
  \sum_{a,b=0}^{2}
  \mathop{{\;\,\sum}'}_{|\mu|,|\nu|\le 1}
  Q_{ab}(\pa^{\mu} \phi,\beta\pa^{\nu} \psi),\\
 &N_{KG}=
  \sum_{a,b=0}^{2}
  \mathop{{\;\,\sum}'}_{|\mu|,|\nu|\le 1}
  Q_{ab}(\pa^{\mu} \psi^{\dagger},\beta\pa^{\nu} \psi),
\end{align*}
and $R_{D}$, $R_{KG}$ are smooth functions of 
$((\pa^{\mu} \psi)_{|\mu|\le 2}, (\pa^{\nu} \phi)_{|\nu|\le 2})$ 
which vanish of cubic order at $(0,0)$.
\end{crl}

\begin{rem}
Roughly saying, the above assertion tells us that the right-hand side in 
(\ref{DKG}) are splitted into two parts: 
The first one is the image of the corresponding linear operator, and the 
second one consists of faster decaying terms which can be regarded as harmless 
remainder when $t \gg 1$. By pushing the first part into the left-hand side, 
we can rewrite (\ref{DKG}) as
\begin{align*}
\left\{
\begin{array}{l}
\dirac_M (\psi-\conjdirac_M \Lambda_D)=N_D+R_D,\\
(\Box+m^2)(\phi-\Lambda_{KG})=N_{KG}+R_{KG}.
\end{array}\right.
\end{align*}
This is what we call the normal form transformation.
\end{rem}

\section{A priori estimate} 
The goal of this section is to get some a priori estimate. 
From now on, we consider only the forward Cauchy problem (i.e., $t>0$) 
since the backward problem can be treated in the same way. 
Let $(\psi,\phi)$ be a solution of 
(\ref{DKG})--(\ref{data_DKG}) for $t \in [0,T)$. We define 
\begin{align*}
 E(T)=\sup_{0\le t <T} \Bigl[&
 \langle t \rangle^{-\delta} \bigl(
  \|\psi(t)\|_{s}+\|\pa \psi(t)\|_{s}+ \|\phi(t)\|_{s}+\|\pa \phi(t)\|_{s}
  \bigr)\\
 &+
\|\psi(t)\|_{s-2}+\|\pa \psi(t)\|_{s-2}
+
\|\phi(t)\|_{s-2}+\|\pa \phi(t)\|_{s-2}\\
&+\sup_{x \in \R^2} \bigl\{\langle t +|x|\rangle 
(|\psi(t,x)|_{s-8}+|\phi(t,x)|_{s-8})\bigr\}
\Bigr],
\end{align*}
where $s \ge 18$ and $0<\delta<1$. Then we have the following.

\begin{prp} \label{prp_apriori}
Let $m \ne 2M$. 
Assume that \eqref{DataSize} is satisfied. Suppose that $E(T)\le 1$. 
There exists a positive constant $C_0$, which is independent of $\eps$ 
and $T$, such that 
\begin{align}
  E(T) \leq C_0(\eps + E(T)^2). 
  \label{mainest}
\end{align}
\end{prp}

We omit the proof of this proposition because it is exactly the same as 
that of the previous works (\cite{katayama}, \cite{sunagawa}, \cite{k_o_s}, 
etc.). The point is that Corollary \ref{cor_normalform} and the commutation 
relation (\ref{commutation}) imply 
\begin{align*}
\left\{
\begin{array}{l}
(\Box+M^2)Z^{\nu}(\psi-\conjdirac_M \Lambda_D)=Z^{\nu}\conjdirac_M(N_D+R_D),\\
(\Box+m^2)Z^{\nu}(\phi-\Lambda_{KG})=Z^{\nu}(N_{KG}+R_{KG})
\end{array}\right.
\end{align*}
with 
\begin{align*}
 |Z^{\nu}\Lambda_*(t,x)| 
 &\lesssim  
 |u|_{[|\nu|/2]+1} (|u|_{|\nu|}+|\pa u|_{|\nu|}),\\
 |Z^{\nu} R_*(t,x)|
   &\lesssim  | u |_{[|\nu|/2]+2}^2 
              \bigl( |u|_{|\nu|+1} + |\pa u|_{|\nu|+1} \bigr),
\end{align*}
and
\begin{align*}
 |Z^{\nu} N_*(t,x)|
  \lesssim 
     \frac{1}{\langle t+|x| \rangle}\bigl|u \bigr|_{[|\nu|/2]+2} 
       \bigl( |u|_{|\nu|+1} + |\pa u|_{|\nu|+1} \bigr),
\end{align*}
where $u=(\psi,\phi)$, and $*$ stands for ``$D$" or ``$KG$". 
Remark that the restriction $s\ge 18$ comes from the relation 
$[(s+1)-2]/2+2 \le s-8$.

\section{Proof of Theorem \ref{thm_main}} 
Now we are ready to prove Theorem \ref{thm_main}. 
First we examine the global existence part of the theorem. 
The inequality (\ref{mainest}) implies 
that there exists a constant $\rho>0$, 
which does not depend on $T$, such that 
$$E(T) \le \rho$$ 
if we choose $\eps$ sufficiently small. 
The unique global existence of the solution for (\ref{DKG})--(\ref{data_DKG}) 
is an immediate consequence of this a priori bound 
and the classical local existence theorem. 

Next we turn to the proof of the existence of the scattering state. 
Remember that 
$$
 \dirac_M (\psi-\conjdirac_M\Lambda_D)=N_D+R_D
$$
with 
\begin{align*}
 &\|(N_D+R_D)(t,\cdot)\|_{H^{s-1}}
 \lesssim 
 \langle t\rangle^{-2+\delta},\\
 &\|\conjdirac_M \Lambda_D(t,\cdot)\|_{H^{s-1}}
 \lesssim 
 \langle t\rangle^{-1+\delta}.
\end{align*}
Now we set
$$
 \psi_0^+:= \psi_0-(\conjdirac_M \Lambda_D)|_{t=0} + 
 \int_{0}^{\infty} U_D(-\tau) (N_D+R_D)(\tau) d\tau
$$
and $\psi^+(t)=U_D(t)\psi_0^+$, 
where $U_D(t)=\exp(-t(\alpha\cdot \nabla_x +iM\beta))$. 
Since the Duhamel formula yields 
\begin{align*} 
 \psi(t)-\conjdirac_M \Lambda_D(t)
 &= 
 U_D(t)(\psi_0 -(\conjdirac_M \Lambda_D)|_{t=0})
 + 
 \int_{0}^{t} U_D(t-\tau) (N_D+R_D)(\tau) d\tau\\
 &=
 \psi^{+}(t) 
 -
 \int_{t}^{\infty} U_D(t-\tau) (N_D+R_D)(\tau) d\tau,
\end{align*}
we have
\begin{align*} 
 \|\psi(t)-\psi^{+}(t)\|_{H^{s-1}}
 &\le  
 \|\conjdirac_M \Lambda_D(t)\|_{H^{s-1}}
 +
 \int_{t}^{\infty} \|(N_D+R_D)(\tau)\|_{H^{s-1}} d\tau\\
 &\lesssim 
 \langle t \rangle^{-1+\delta} 
 +
 \int_{t}^{\infty}  \langle \tau \rangle^{-2+\delta}  d\tau\\
 &\lesssim 
 \langle t \rangle^{-1+\delta}. 
\end{align*}
As for the Klein-Gordon component, we just have to set 
 \begin{align*} 
 &\phi_0^+
 =
 \phi_0-\Lambda_{KG}\bigr|_{t=0} 
 +
 \int_{0}^{\infty} 
 \frac{\sin \left(-\tau \OOmega_m \right)}{\OOmega_m }
 (N_{KG}+R_{KG})(\tau,\cdot) d\tau,\\
 &\phi_1^+
 =
 \phi_1-\pa_t \Lambda_{KG} \bigr|_{t=0} 
 +
 \int_{0}^{\infty} 
 \left(\cos(-\tau \OOmega_m )\right) (N_{KG}+R_{KG})(\tau,\cdot) d\tau
\end{align*}
with $\OOmega_m=(m^2-\Delta)^{1/2}$. 
\qed


\end{document}